\documentclass[onecolumn]{elsarticle}

\journal{Journal of Mathematical Analysis and Applications}

\usepackage[english]{babel}
\usepackage{amsfonts}
\usepackage{amsmath}
\usepackage{amscd}
\usepackage[active]{srcltx} 
\usepackage{latexsym}
\usepackage{amssymb}
\usepackage[latin1]{inputenc}
\usepackage{enumerate}

\newtheorem{thm}{Theorem}[section]
 \newtheorem{cor}[thm]{Corollary}
 \newtheorem{lema}[thm]{Lemma}
 \newtheorem{prop}[thm]{Proposition}
 \newdefinition{Def}{Definition}
 \newtheorem{obs}[thm]{Remark}
 \newproof{Dem}{Proof}

 \numberwithin{equation}{section}

\newcommand{\PI}[2]{\left\langle \,#1 , #2\, \right\rangle}
\newcommand{\PIW}[2]{\left\langle \,#1 , #2\, \right\rangle_W}

\newcommand{\NW}[1]{\Vert #1 \Vert_{W}^{2}}
\newcommand{\NC}[1]{\Vert #1 \Vert}

\newcommand{\NWW}[1]{\Vert #1 \Vert_{W}}

\newcommand{\Nphis}[1]{\Vert #1 \Vert_{p}}
\newcommand{\Nphiss}[1]{\Vert #1 \Vert_{p}}
\newcommand{\Nphissw}[1]{\Vert #1 \Vert_{p,W}}

\newcommand{\ra}{\rightarrow}

\newcommand{\St}{\mathcal{S}}
\newcommand{\HH}{\mathcal{H}}
\newcommand{\M}{\mathcal{M}}
\newcommand{\N}{\mathcal{N}}

\newcommand{\CC}{\mathbb{C}}
\newcommand{\Q}{\mathcal{Q}}

\newcommand{\ol}{\overline}

\DeclareMathOperator{\rank}{rank}

\DeclareMathOperator{\argmin}{argmin}

\begin{document}

\begin{frontmatter}

\title{Weighted Procrustes problems}

\author[FI]{Maximiliano Contino}
\ead{mcontino@fi.uba.ar}

\author[FI,IAM]{Juan Ignacio Giribet\corref{ca}}
\ead{jgiribet@fi.uba.ar}

\author[FI,IAM]{Alejandra Maestripieri}
\ead{amaestri@fi.uba.ar}

\cortext[ca]{Corresponding author}
\address[FI]{%
Facultad de Ingenier\'{\i}a, Universidad de Buenos Aires\\
Paseo Col\'on 850 \\
(1063) Buenos Aires,
Argentina 
}

\address[IAM]{%
Instituto Argentino de Matem\'atica ``Alberto P. Calder\'on'' \\ CONICET\\
Saavedra 15, Piso 3\\
(1083) Buenos Aires, 
Argentina }

\begin{abstract}
Let $\HH$ be a Hilbert space, $L(\HH)$ the algebra of bounded linear operators on $\HH$ and $W \in L(\HH)$ a positive operator such that $W^{1/2}$ is in the p-Schatten class, for some $1 \leq p< \infty.$ Given $A \in L(\HH)$ with closed range and $B \in L(\HH),$ we study the following weighted approximation problem: analize the existence of 
$$\underset{X \in L(\HH)}{min}\Nphissw{AX-B},$$
where $\Nphissw{X}=\Nphiss{W^{1/2}X}.$

In this paper we prove that the existence of this minimum is equivalent to a compatibility condition between $R(B)$ and $R(A)$ involving the weight $W,$ and we characterize the operators which minimize this problem as $W$-inverses of $A$ in $R(B).$
\end{abstract}


\end{frontmatter}

\section{Introduction}
One problem of interest in Signal and Image Processing is to find low dimensional models that approximate, in some sense, given data \cite{Dai2010, Erikssonand2010}. In particular many of these problems can be posed as follows:
given a matrix $B\in \CC^{n\times n}$, with $\rank(B)\geq k$, for $k\in\mathbb{N}$ satisfying $k<n$, find a matrix $Y_0\in \CC^{n\times n}$ with $\rank(Y_0)=k$ such that,
$$
Y_0=\underset{\{Y\in\CC^{n\times n}: \ \rank(Y)=k\}}{\argmin} f(Y-B),
$$
for some {\it cost} function $f:\CC^{n\times n}\rightarrow\mathbb{R}$. Due to its intractability, usually this problem is studied by relaxating the constraint on the rank of $Y,$ which under certain conditions, turns out to be an exact relaxation. For this, the  factorization $Y=AX$ is used, with $A\in\CC^{n\times k}$, $X\in\CC^{k\times n}$. Assume that the cost function is given by the Frobenius norm $\Vert \cdot \Vert_F$, now we are interested in the following problem:
\begin{eqnarray}\label{prob_rango}
Y_0=\underset{\{X\in\CC^{k\times n}, \ A\in\CC^{n\times k}\}}{\argmin} \|AX-B\|_F.
\end{eqnarray}
In fact, suppose that $A_0\in\CC^{n\times k}, X_0\in\CC^{k\times n}$ satisfy
\begin{eqnarray} \label{A0X0-B}
\|A_0X_0-B\|_F=\min_{\{X\in\CC^{k\times n}, \ A\in\CC^{n\times k}\}} \|AX-B\|_F,
\end{eqnarray}
then, 
\begin{eqnarray} \label{A0X-B}
\|A_0X_0-B\|_F=\min_{X\in\CC^{k\times n}} \|A_0X-B\|_F.
\end{eqnarray}

If a (positive) weight is introduced in equation \eqref{A0X-B}, or if the Frobenius norm is replaced by another unitary invariant norm, the same problem can be studied. 

This work is devoted to study an extension of problem \eqref{A0X-B} in abstract Hilbert spaces. More specifically, we study the following approximation problem. Given $A\in L(\HH)$ with closed range, $B\in L(\HH)$ and $W\in L(\HH)$ a positive operator, we analyze the conditions for the existence of
\begin{equation}\label{inf}
\min_{X \in L(\HH)}\|W^{1/2}(AX-B)\|_p,
\end{equation}
for $1\leq p <\infty$, where $\Nphiss{\cdot}$ is the p-Schatten norm.

There are several examples of these minimization problems in Control Theory and Signal Processing \cite{Dejoudi, Singler}. Similar problems also arise in Quantum Chemistry, for example in the orthogonalization process of L\"owdin \cite{Aiken2, Low1}, or in the approximation of the Hamiltonian operator \cite{Gold1, Gold2, Mah1}.

The existence of minimum of $\NC{AX-B}_{p}$ in Hilbert spaces, was studied in \cite{Mah} using differentiation techniques and also in \cite{Nashed}, where a connection between $p$-Schatten norms and the order in $L(\HH)^+$ (the cone of semidefinite positve operators) is established. However, the introduction of a weight $W\in L(\HH)^+$ plays an important role, since we are introducing on $\HH$ a semi-inner product associated to $W$ for which $\HH$ is no longer a Hilbert space, unless $W$ is invertible. In this case, the existence of a suitable orthogonal projection is not guaranteed. In fact the existence of a $W$-orthogonal projection onto $R(A)$ depends on the relationship between the weight $W$ and the closed subspace $R(A)$.

The notion of compatibility, defined in \cite{CMSSzeged} and developed later in \cite{WGI, Shorted3, Spline}, has its origin in the work of Z. Pasternak-Winiarski \cite{Pasternak}. In that work the author studied, for a fixed subspace $\St$, the analiticity of the map $W \to P_{W,\St}$ which associates to each positive invertible operator $W$ the orthogonal projection onto $\St$ under the (equivalent) inner product $\PI{x}{y}_W=\PI{Wx}{y}$, for $x, y \in \HH$. The notion of compatibility appears when $W$ is allowed to be any positive semidefinite operator, not necessarily invertible (and even, a selfadjoint bounded linear operator). More precisely, $W$ and $\St$ are said to be $\it  compatible $ if there exists a (bounded linear) projection $Q$ with range $\St$ which satisfies $WQ=Q^*W.$ If $W$ is positive and invertible or $\HH$ has finite dimension, there exists a unique projection onto $\St$ which is $W$- selfadjoint \cite{CMSSzeged}. In general, it may happen that there is no such $Q$ or that there is an infinite number of them.
However, there exists an angle condition between $\St^\bot$ and $\ol{W(\St)}$ which determines the existence of these projections \cite{Dixmier}. In fact, the existence of such projections is related with the existence of minimum of equation \eqref{inf}.

\vspace{0,3cm}
The contents of the paper are the following. In section 2, some characterizations of the compatibility of the pair $(W,R(A))$ are given. Also some properties of shorted operators and  compressions and its connection with compatibility is stated. Finally, the concept of $W$-inverses of an operator $A$ in the range of an operator $B$, and some properties are presented. 

For the sake of simplicity, in section 3, we study problem \ref{inf} when $B=I.$ We prove that the infimum of the set $\{(AX-I)^{*}W(AX-I): X \in L(\HH)\}$ (where the order is the one induced by the cone of positive operators), always exists and is equal to $W_{/R(A)},$ the shorted operator of $W$ to $R(A)$. We also prove that the existence of the minimum of the previous set is equivalent to the compatibility of the pair $(W, R(A)).$  We characterize the operators which minimize this problem, which are the $W$-inverses of $A.$ Finally, it is shown that if $W^{1/2}$ is in the $p$-Schatten class, for some $1 \leq p< \infty,$ the existence of the minimum of the set $\{ \Nphiss{W^{1/2}(AX-I)} : X \in L(\HH)\}$ is also equivalent to the compatibility of the pair $(W, R(A))$. In this case,  set of solutions of \ref{inf} are the $W$-inverses of $A.$

In section 4, we prove similar results for an arbitrary operator $B \in L(\HH)$, where the existence of the minimum of the set $\{ \Nphiss{W^{1/2}(AX-B)} : X \in L(\HH)\},$ is equivalent to the compatibility condition $R(B) \subseteq R(A)+W(R(A))^{\perp}.$ In this case, the minimizers are the $W$-inverses of $A$ in $R(B)$. 

\section{Preliminaries}
In the following $\HH$ denotes a separable complex Hilbert space, $L(\HH)$ is the algebra of bounded linear operators from $\HH$ to $\HH$, and $L(\HH)^{+}$ the cone of semidefinite positive operators. $GL(\HH)$ is the group of invertible operators in $L(\HH),$ $CR(\HH)$ is the subset of $L(\HH)$ of all operators with closed range.
For any $A \in L(\HH),$ its range and nullspace are denoted by $R(A)$ and $N(A)$, respectively. Finally,  $A^{\dagger}$ denotes the  Moore-Penrose inverse of the operator $A \in L(\HH).$

Given two closed subspaces $\M$ and $\N$ of $\HH,$ $\M \dot{+} \N$ denotes the direct sum of $\M$ and $\N$. 
If $\HH$ is decomposed as a direct sum of closed subspaces $\HH=\M \dot{+} \N,$ the projection onto $\M$ with nullspace $\N$ is denoted by $P_{\M {\mathbin{\!/\mkern-3mu/\!}} \N},$ and $P_{\M} = P_{\M {\mathbin{\!/\mkern-3mu/\!}} \M^{\perp}}.$ Also, $\Q$ denotes the subset of $L(\HH)$ of oblique projections, i.e. $\Q=\{Q \in L(\HH): Q^{2}=Q\}.$ 

Given $W\in L(\HH)^{+}$ and a closed subspace $\St$ of $\HH,$ the pair $(W,\St)$ is $\it{compatible}$ if there exists $Q\in \Q$ with $R(Q)=\St$ such that $WQ=Q^{*}W.$ The last condition means that $Q$ is $W$-hermitian, in the sense that $\PIW{Qx}{y}=\PIW{x}{Qx},$ for every $x, y \in \HH,$ where $\PIW{x}{y}=\PI{Wx}{y}$ defines a semi-inner product on $\HH.$ 

The $W$-orthogonal complement of $\St$ is $$\St^{\perp_{W}}=\{x \in \HH: \PI{Wx}{y}=0, \ y\in \St\}=W^{-1}(\St^{\perp}).$$

The next theorem, proven in \cite[Prop.~3.3]{CMSSzeged}, allows us to characterize the compatibility of the pair $(W,\St).$

\begin{thm} \label{Comp 1} Given $W\in L(\HH)^{+}$ and a closed subspace $\St \subseteq \HH,$  the pair $(W,\St)$ is compatible if and only if $$\HH=\St+\St^{\perp_{W}}.$$\end{thm}

\vspace{0,3cm}

 We now give the definitions of $W$-least squares solution of the equation $Az=x.$

\begin{Def} Given $A\in CR(\HH),$ $W\in L(\HH)^{+}$ and $x\in \HH,$ $u\in \HH$ is a $W$-least squares solution or $W$-$LSS$ of $Az=x,$ if $$\NWW{Au-x}\leq\NWW{Az-x}, \mbox{ for every } z \in \HH,$$ where $\NW{x}=\PI{Wx}{x}$ is the seminorm associated to $W$. \end{Def}

\begin{thm}
\label{thmWLSS} Given $A\in CR(\HH),$ $W\in L(\HH)^{+}$ and  $x \in \HH,$ then the following statements hold:
\begin{enumerate}
\item [i)] There exists a $W$-$LSS$ of $Az =x$ if and only if $x \in R(A)+R(A)^{\perp_W},$ 
\item [ii)] $u_0$ is a $W$-$LSS$ of $Az=x$ if and only if $$A^{*}W(Au_0-x)=0.$$ 
\end{enumerate}
\end{thm}

\begin{Dem} For item i) see \cite{Spline}, and for item ii) see \cite[Remark~5.2]{WGI}.\end{Dem}

The following is a well known result due to R. Douglas \cite{Douglas} about range inclusion and factorizations of operators. In the following we use the operator order induced by $L(\HH)^+$, i.e., $A \leq B$ if $B-A \in L(\HH)^+$.
\begin{thm} \label{Teo8} Let $Y, Z \in L(\HH)$, the following conditions are equivalent:
\begin{enumerate}
\item [i)] $R(Z)\subseteq R(Y),$
\item [ii)] there exists a positive number $\lambda$ such that $ZZ^{*}\leq \lambda YY^{*},$
\item [iii)] there exists $D\in L(\HH)$ such that $Z=YD.$
\end{enumerate}
\end{thm}

In this case there exists a unique solution $D_0$ of the equation $Z=YX,$ such that $R(D_0) \subseteq N(Y)^{\perp}.$ Moreover $\Vert D_0 \Vert = inf \{\lambda: ZZ^{*}\leq \lambda YY^{*} \}.$\\

\vspace{0,3cm}

 In \cite{Mitra} S. K. Mitra and C. R. Rao introduced the notion of the $W$-inverse of a matrix. We extend the definition in the following way. 

\begin{Def} Given $A \in CR(\HH),$ $B \in L(\HH)$ and $W\in L(\HH)^{+},$ $X_0 \in L(\HH)$ is a $W$-inverse of $A$ in $R(B),$ if for each $x \in \HH$, $X_0x$ is a $W$-$LSS$ of $Az=Bx,$ i.e. $$\NWW{AX_0x-Bx}\leq\NWW{Az-Bx}, \mbox{ for every } x, z \in \HH.$$ \end{Def} 

When $B=I,$ $X_0$ is called a $W$-inverse of $A$, see \cite{WGI}. The next theorem shows that there is a close relationship between $W$-inverses and $W$-LSS solutions. 

\begin{thm}
\label{thmWinversa} Given $A\in CR(\HH), B \in L(\HH)$ and $W\in L(\HH)^{+},$ the following conditions are equivalent:
\begin{enumerate} 
\item [i)] The operator $A$ admits a $W$-inverse in $R(B),$
\item [ii)] $R(B) \subseteq R(A) + R(A)^{\perp_{W}},$ 
\item [iii)] the normal equation $A^{*}WAX=A^{*}WB$ admits a solution. 
\end{enumerate}
\end{thm}

\begin{Dem}
$i) \Leftrightarrow iii):$ If $X_0$ is a $W$-inverse of $A$ in $R(B)$ then,
$$\NWW {AX_0x-Bx} \leq \NWW {Az-Bx}, \mbox{ for every } x, z \in L(\HH).$$
Or equivalently, $X_0x$ is a $W$-$LSS$ of $Az=Bx,$ for every $x \in \HH.$ Or, by Theorem \ref{thmWLSS}, 
$$A^{*}W(AX_0-B)x=0, \mbox{ for every } x \in \HH,$$
so that $X_0$ is a solution of the normal equation.
The converse follows in a similar way, applying Theorem \ref{thmWLSS}.

$ii) \Leftrightarrow iii):$ If $R(B) \subseteq R(A) + R(A)^{\perp_{W}},$ applying $A^{*}W$ to both sides of the inclusion, 
$$R(A^{*}WB) \subseteq R(A^{*}WA).$$
Then by Theorem \ref{Teo8}, the normal equation admits a solution.
The converse follows easily.
\qed
\end{Dem}

\begin{cor}\label{Winversa 3} 
If $R(B) \subseteq R(A) + R(A)^{\perp_{W}},$ then the set of $W$-inverses of $A$ in $R(B)$ is the set of solutions of the equation $A^{*}WAX=A^{*}WB$, or equivalently the affine manifold
$$(A^{*}WA)^{\dagger}A^{*}WB+ \{L \in L(\HH): R(L) \subseteq N(A^{*}WA)\}.$$
\end{cor}
\vspace{0,3cm}

 Given a positive operator $W\in L(\HH)^{+}$ and a closed subspace $\St \subseteq \HH$ the notion of shorted operator of $W$ to $\St,$ was introduced by M. G. Krein in \cite{Krein} and later rediscovered by W. N. Anderson and G. E. Trapp who proved in \cite{Shorted2}, that the set
$\{ X \in L(\HH): \ 0\leq X\leq W \mbox{ and } R(X)\subseteq \St^{\perp}\}$ has a maximum element.

\begin{Def} The shorted operator of $W$ to $\St$ is defined by
$$W_{/\St}=\mbox{max } \{ X \in L(\HH): \ 0\leq X\leq W \mbox{ and } R(X)\subseteq \St^{\perp}\}.$$

The $\St$-compression $W_{\St}$ of $W$ is the (positive) operator defined by $$W_{\St}=W-W_{/\St}.$$\end{Def}

For many results on the notions of shorted operators, the reader is referred to \cite{Shorted1} and \cite{Shorted2}. 
Next we collect some results regarding $W_{/\St}$ and $W_{\St}$ which will be used in the rest of this work.

\begin{thm} \label{TeoShorted}
Let $W\in L(\HH)^{+}$ and $\St \subseteq \HH$ a closed subspace. Then  
\begin{enumerate} 
\item [i)] $W_{/\St}=\mbox{ inf } \{ E^{*}WE: E^2=E, \ N(E)=\St\};$ in general, this infimum is not attained,
\item [ii)] $R(W) \cap \St^{\perp} \subseteq R(W_{/\St}) \subseteq R(W^{1/2}) \cap \St^{\perp}$, and  $\overline{N(W)+\St} \subseteq N(W_{/\St})=W^{-1/2}(\overline{W^{1/2}(\St)}),$
\item [iii)] $N(W_{\St})=W^{-1}(\St^{\perp})$ and $W(\St) \subseteq R(W_{\St}) \subseteq \overline{W(\St)}.$
\end{enumerate}
\end{thm}

 The reader is referred to \cite{Shorted2} and \cite{Shorted3} for the proof of these facts.
In \cite{Shorted3} the next results were stated.

\begin{thm} \label{TeoShorted2} Let $W\in L(\HH)^{+}$ and $\St \subseteq \HH$ be a closed subspace. The following conditions are equivalent:
\begin{enumerate} 
\item [i)] The pair $(W,R(A))$ is compatible,
\item [ii)] $W_{/\St}=\mbox{min } \{ E^{*}WE: E^2=E, \ N(E)=\St\},$
\item [iii)] $R(W_{/\St})=R(W) \cap \St^{\perp} \mbox{ and } N(W_{/\St})=N(W)+\St.$
\end{enumerate}
\end{thm}

\vspace{0,3cm}
\begin{Def} Let $T\in L(\HH)$ be a compact operator. By $\{\lambda_k(T)\}_{k\geq1}$ we denote the eigenvalues of $\vert T \vert = (T^{*}T)^{1/2},$ where each eigenvalue is repeated according to its multiplicity. Let $1\leq p < \infty,$ we say that $T$ belongs to the p-Schatten class $S_p,$ 
if $$\sum_{k\geq1}^{} \lambda_{k}(T) ^{p}<\infty,$$ and the $p$-Schatten norm is given by 
$$\Nphiss{T}= (\sum_{k\geq1}^{} \lambda_{k}(T) ^{p})^{1/p},$$
\end{Def}

The reader is referred to \cite{Ringrose, Simon} for a detailed exposition of these topics. The Schatten norms are unitary invariant norms. More generally, 

\begin{Def} A norm $\vert \vert \vert \cdot \vert \vert \vert$ on a non-zero ideal $\mathcal{J}$ of $L(\HH)$ is called unitarily invariant if 
$$\vert \vert \vert UTV \vert \vert \vert= \vert \vert \vert T \vert \vert \vert,$$ for any unitary operators $U,V \in L(\HH)$ and $T \in \mathcal{J}$. \end{Def}

\begin{lema} \label{lemaNUI} Every unitarily invariant norm $\vert \vert \vert \cdot \vert \vert \vert$  on a non-zero ideal $\mathcal{J}$ of $L(\HH)$ is symmetric, i.e., 
$$\vert \vert \vert S_1 T S_2 \vert \vert \vert \leq  \Vert S_1 \Vert  \ \vert \vert \vert T  \vert \vert \vert \ \Vert S_2 \Vert,$$ 
for every $S_1, S_2 \in L(\HH)$ and $T \in \mathcal{J}.$\end{lema}

\begin{Dem} See \cite[Lemma.~2.1]{Kitta}.\end{Dem}

The following result will be useful to study problem \ref{inf}. A more general result can be found in \cite[Proposition 2.5]{Nashed}.

\begin{prop}\label{Prop Nashed} Let $\vert \vert \vert \cdot \vert \vert \vert$ be an unitarily invariant norm on a non-zero ideal $\mathcal{J}$ of $L(\HH),$ and $S, T \in \mathcal{J}.$ Then, $$\mbox{ if }T^{*}T\leq S^{*}S \mbox{ then } \vert \vert \vert T \vert \vert \vert \leq \vert \vert \vert S \vert \vert \vert.$$ \end{prop}
 
\begin{Dem} 
If $T^{*}T \leq S^{*}S,$ by Theorem \ref{Teo8}, there exists an operator $R$ with $\Vert R \Vert \leq 1$ such that $T^{*}=S^{*}R,$ then using Lemma \ref{lemaNUI} we have 
$\vert \vert \vert T \vert \vert \vert=\vert \vert \vert T^{*} \vert \vert \vert = \vert \vert \vert S^{*}R \vert \vert \vert  \leq \vert \vert \vert S^{*} \vert \vert \vert \ \Vert R \Vert \leq \vert \vert \vert S^{*} \vert \vert \vert=\vert \vert \vert S \vert \vert \vert.$ \hfill 
\qed
\end{Dem}

\vspace{0,5cm}
Finally, we give a definition of a derivative that will be instrumental to prove some results stated in Section 3. 

\begin{Def} Let $(\mathcal{E}, \NC{\cdot})$ be a Banach space and $f: \mathcal{E} \rightarrow \mathbb{R}.$ Let $\phi \in [0, 2\pi)$ and $h>0$, then the $\phi-$directional derivative of $f$ at a point $x\in \mathcal{E}$ in direction $y \in \mathcal{E}$ is defined by
$$ D_{\phi}f(x,y)=lim_{h\rightarrow 0^{+}} \frac{f(x+h e^{i\phi}y)-f(x)}{h}.$$\end{Def}

\begin{thm} \label{TeoD} Let $G_p: S_p \rightarrow \mathbb{R}^{+},$ $1\leq p<\infty,$ $G_p(X)=\Nphis{X}^{p},$ and let $X, Y \in S_p.$ Then, 
\begin{itemize}
\item [i)] For $1< p<\infty,$ $G_p$ has a $\phi-directional$ derivative given by 
$$D_{\phi}G_p(X,Y)=p \ Re \ [e^{i\phi} tr ( \vert X \vert^{p-1} U^{*}Y)],$$ for all $\phi \in [0, 2\pi).$
\item[ii)]  For $p=1,$ $G_1$ has a $\phi-directional$ derivative given by 
$$D_{\phi}G_1(X,Y)=Re \ [e^{i\phi} tr(U^{*}Y)] + \Vert P_{N(X^{*})}YP_{N(X)}\Vert_{1},$$
for all $\phi \in [0, 2\pi),$
\end{itemize}
where $Re (z)$ is the real part of a complex number $z$, $tr(T)$ denotes the trace of the operator $T$ and $X=U\vert X\vert,$ is the polar descomposition of the operator $X,$ with $U$ the partial isometry such that $N(U)=N(X).$\end{thm}
\begin{Dem} See \cite[Theorem 2.1]{Aiken} and \cite[Theorem 2.1]{Drago}. \end{Dem}

\begin{lema} \label{LemaM2} Let $(\mathcal{E}, \NC{\cdot})$ be a Banach space and $f: \mathcal{E} \rightarrow \mathbb{R},$ such that $f$ has a $\phi-directional$ derivative for every $\phi \in [0, 2\pi),$ at every point $x \in \mathcal{E}$ and in every direction $y \in \mathcal{E}$.  
If $f$ has a global minimum at $x_0 \in \mathcal{E},$ then 
$$\underset{0 \leq \phi < 2\pi}{inf} (D_{\phi}f(x_0,y)) \geq 0, \mbox{ for every } y \in \mathcal{E}.$$\end{lema}

\begin{Dem} See \cite[Theorem 2.1]{Mecheri}.\end{Dem}

\section{Weighted least squares problems}

Given $W \in L(\HH)^{+}$ such that $W^{1/2} \in S_p$ for some $p$ with $1 \leq p < \infty,$ consider the operator seminorm associated to $W$, $$\Nphissw{X}=\Nphiss{W^{1/2}X},$$
for $X\in L(\HH).$
We study the following approximation problem: given $A\in CR(\HH)$ and $B\in L(\HH)$, analize the existence of
$$\underset{X \in L(\HH)}{min}\Nphissw{AX-B}.$$

In this section we study the case when $B=I$, i.e., we study the problem
\begin{equation}
\underset{X \in L(\HH)}{min} \Nphissw{AX-I}. \label{eq1}
\end{equation}

To study problem \eqref{eq1} we introduce the following associated problem: given $W \in L(\HH)^{+}$ and $A\in CR(\HH)$, define $F:L(\HH) \ra L(\HH)^+,$ $$F(X)=(AX-I)^{*}W(AX-I),$$
and analize the existence of
\begin{equation}
\underset{X \in L(\HH)}{inf} F(X), \label{eq2}
\end{equation}
in the order induced in $L(\HH)$ by the cone of positive operators.
The next result shows that the infimum of equation \eqref{eq2} always exists and coincides with the shorted operator of $W$ to $R(A).$

\begin{prop} \label{Teo7} Let $A \in CR(\HH)$ and $W \in L(\HH)^{+},$ then the infimum of problem \eqref{eq2} exists and
$$\underset{X \in L(\HH)}{inf}F(X)=W_{/R(A)}.$$\end{prop}

\begin{Dem} 
Let $X\in L(\HH), $ writing $W=W_{/R(A)}+W_{R(A)},$ it follows that
$$(I-AX)^{*}W(I-AX)=W_{/R(A)}+(I-AX)^{*}W_{R(A)}(I-AX)\geq W_{/R(A)},$$
because $R(A) \subseteq N(W_{/R(A)})$ (see Theorem \ref{TeoShorted}) and then $W_{/R(A)}(I-AX)=W_{/R(A)}=(I-AX)^{*}W_{/R(A)}.$
Hence $W_{/R(A)}$ is a lower bound of $F(X).$

If $C\geq 0$ is any other lower bound of $F(X),$ then $$C\leq F(X),\mbox { for every } X \in L(\HH).$$
In particular, $$C\leq E^{*}WE,$$ 
where $E$ is any projection such that $N(E)=R(A)$. In fact $R(I-E)= N(E)=R(A),$ then by Theorem \ref{Teo8}, there exists $X_0 \in L(\HH),$ such that $(I-E)=AX_0,$ i.e. $(-E)=AX_0-I.$
Therefore, by Theorem \ref{TeoShorted} $$C \leq inf \{E^{*}WE : \ E^{2}=E, \ N(E)=R(A) \}= W_{/R(A)}.$$ 
Thus, $$W_{/R(A)}=\underset{X \in L(\HH)}{inf}F(X).$$ \hfill \qed
\end{Dem}

\begin{thm}  \label{Teo1}Let $A \in CR(\HH)$ and $W \in L(\HH)^{+}.$ Then
problem \eqref{eq2} has a minimum, i.e., there exists $X_0\in L(\HH)$ such that 
$$F(X_0)=\underset{X \in L(\HH)}{min}F(X)=W_{/R(A)}$$
if and only if the pair $(W, R(A))$ is compatible.\end{thm}

\begin{Dem}
If problem \eqref{eq2} has a minimum, from Proposition \ref{Teo7} it holds that there exists $X_0\in L(\HH)$ such that
$$F(X_0)=\underset{X \in L(\HH)}{min}F(X)=W_{/R(A)}.$$ 
Writing again $W=W_{/R(A)}+W_{R(A)}$, it follows that
$$W_{/R(A)}=F(X_0)=W_{/R(A)}+(AX_0-I)^{*}W_{R(A)}(AX_0-I).$$
Therefore $$(AX_0-I)^{*}W_{R(A)}(AX_0-I)=0,$$ then,
$$W_{R(A)}^{1/2}(AX_0-I)=0,$$ and then by Theorem \ref{TeoShorted}
$$R(AX_0-I) \subseteq N(W_{R(A)})=W^{-1}(R(A)^{\perp}).$$
Therefore $$W(R(AX_0-I))\subseteq R(A)^{\perp} \cap R(W).$$
Then $$R(W_{/R(A)})=R(F(X_0))\subseteq (AX_0-I)^{*} (R(A)^{\perp} \cap R(W))=R(A)^{\perp} \cap R(W),$$
because $A^{*}(R(A)^{\perp})=0.$ Then, $R(W_{/R(A)})= R(A)^{\perp} \cap R(W),$ because $R(A)^{\perp} \cap R(W)$ is always contained in $R(W_{/R(A)})$ (see Theorem \ref{TeoShorted}).

Also $x \in N(W_{/R(A)})$ if and only if $W^{1/2}(AX_0-I)x=0,$ or equivalently $(AX_0-I)x \in N(W).$
In this case $x \in N(W)+R(A),$ and then $$N(W_{/R(A)})=N(W)+R(A),$$
because $N(W)+R(A)$ is always contained in $N(W_{/R(A)})$ (see Theorem \ref{TeoShorted}). Therefore 
$R(W_{/R(A)}) = R(A)^{\perp} \cap R(W)$ and $N(W_{/R(A)})=N(W)+R(A)$ and by Theorem \ref{TeoShorted2}, the pair $(W, R(A))$ is compatible.

Conversely, if the pair $(W, R(A))$ is compatible, then by Theorem \ref{TeoShorted2}, $$W_{/R(A)}=min \{E^{*}WE: \ E^2=E, \ N(E)=R(A)\}.$$ Let $E_0$ be such that $E_0^{2}=E_0,$ $N(E_0)=R(A)$ and $W_{/R(A)}=E_0^{*}WE_0.$ Consider $X_0=A^{\dagger}(E_0-I),$ then $E_0=AX_0-I$ and $F(X_0)=W_{/R(A)}.$ 
\hfill \qed
\end{Dem}

\vspace{0,3cm} If the pair $(W, R(A))$ is compatible then by Theorem \ref{Teo1}, problem \eqref{eq2} attains a minimum, i.e., there exists $U_0 \in L(\HH)$ such that $F(U_0)=W_{/R(A)}.$ Consider the set 
$$M=\{X\in L(\HH): F(X)=W_{/R(A)}\}.$$ The next proposition gives a characterization of the elements of $M.$

\begin{prop} \label{Teo2} Let $A \in CR(\HH)$ and $W \in L(\HH)^{+}$ such that the pair $(W, R(A))$ is compatible. The following conditions are equivalent:
\begin{itemize}[leftmargin=*]
\item [i)]  $X_0 \in M,$ i.e. $F(X_0)=\underset{X \in L(\HH)}{min} F(X),$
\item [ii)] $X_0$ is a $W$-inverse of $A,$ 
\item [iii)] $X_0$ is a solution of the normal equation $$A^{*}W(AX-I)=0.$$
\end{itemize}\end{prop} 

\begin{Dem} $i) \Leftrightarrow ii):$ If $X_0$ is such that $F(X_0)\leq F(X)$, for every $X\in L(\HH)$, then 
$$ \NW {AX_0x-x} \leq \NW {AXx-x}, \mbox{ for every } x\in\HH \mbox{ and every } X\in L(\HH).$$  
For every $x \in \HH$, given $z\in \HH,$ let $X\in L(\HH)$ such that  $z=Xx.$ Then $$ \NW {AX_0x-x} \leq \NW {Az-x}, \mbox{ for every } x\in \HH \mbox{ and every } z\in \HH.$$Therefore $X_0$ is a $W$-inverse of $A.$ The converse is similar.


The equivalence $ii) \Leftrightarrow iii)$ was established in Theorem \ref{thmWinversa}, for $B=I.$
\hfill \qed
\end{Dem}

\vspace{0,3cm}
\begin{obs} Let $A \in CR(\HH)$ and $W \in L(\HH)^{+}$ such that $W^{1/2} \in S_p,$ for some $p$ with $1 \leq p < \infty,$ then
$$ \underset{X \in L(\HH)}{inf} \Nphissw {AX-I} \geq \Nphiss {W_{/R(A)}^{1/2}}.$$

In fact, by Proposition \ref{Teo7}, $$(AX-I)^{*}W(AX-I)\geq W_{/R(A)}, \mbox{ for every } X \in L(\HH).$$
By Proposition \ref{Prop Nashed} we get that $$\underset{X \in L(\HH)}{inf} \Nphissw {AX-I} \geq  \Nphiss {W_{/R(A)}^{1/2}}.$$

\end{obs}

 The next result proves the equivalence between the existence of a minimum of problem \eqref{eq1} and the compatibility of the pair $(W,R(A)).$

\vspace{0,3cm}
\begin{thm} \label{Teo3} Let $A \in CR(\HH)$ and $W \in L(\HH)^{+},$ such that $W^{1/2} \in S_p,$ for some $p$ with $1 \leq p < \infty$. Problem \eqref{eq1} has a minimum if and only if the pair $(W,R(A))$ is compatible. 

In this case,
$$\underset{X \in L(\HH)}{min}  \Nphissw {AX-I} = \Nphiss {W_{/R(A)}^{1/2}}.$$
Moreover, $X_0 \in L(\HH)$ satisfies
$$\Nphissw {AX_0-I} = \Nphiss {W_{/R(A)}^{1/2}},$$
if and only if $X_0$ is a $W$-inverse of $A.$
\end{thm}
\begin{Dem} If the pair $(W, R(A))$ is compatible, then by Theorem \ref{Teo1}, there exists $X_0\in L(\HH)$ such that $F(X_0)=min_{X\in L(\HH)}F(X)=W_{/R(A)},$ i.e. $$W_{/R(A)}=F(X_0)\leq F(X), \mbox{ for every } X\in L(\HH).$$ Since $W^{1/2} \in S_p,$  by Proposition \ref{Prop Nashed},  $$\Nphiss {W_{/R(A)}^{1/2}}=\Nphiss{W^{1/2}(AX_0-I)}=\Nphissw{AX_0-I}\leq \Nphissw{AX-I}, \mbox{ for every } X\in L(\HH),$$
then $$\underset{X \in L(\HH)}{min} \Nphissw {AX-I} = \Nphissw {AX_0-I} = \Nphiss {W_{/R(A)}^{1/2}}.$$ 
 
To prove the converse, for $1\leq p < \infty,$ consider  
$F_p: S_p \rightarrow \mathbb{R}^{+},$ $$F_p(X)=\Nphiss{W^{1/2}(AX-I)}^{p}.$$
By Theorem \ref{TeoD}, $F_p$ has a $\phi-directional$ derivative for all $\phi \in [0, 2\pi).$ Then it is easy to check that, for every $\ X, \ Y\in L(\HH)$ and $\phi \in [0, 2\pi),$
$$D_{\phi}F_p (X,Y)=D_{\phi}G_p (W^{1/2}(AX-I),W^{1/2}AY),$$
where $G_p(X)=\Nphis{X}^{p}.$

Suppose that problem \eqref{eq1} admits a minimum, i.e. there exists $X_0 \in L(\HH),$ a global minimum  of $\Nphissw{AX-I}.$ Then $X_0$ is a global minimum of $F_p$ and, by Lemma \ref{LemaM2}, we have
$$\underset{0 \leq \phi < 2\pi}{inf} (D_{\phi}F_p(X_0,Y)) \geq 0, \mbox{ for every } Y \in L(\HH).$$
Let $W^{1/2}(AX_0-I)=U\vert W^{1/2}(AX_0-I)\vert$ be the polar descomposition of the operator $W^{1/2}(AX_0-I),$ with $U$ a partial isometry with $N(U)=N(W^{1/2}(AX_0-I))$ and $R(U)=\overline{R(W^{1/2}(AX_0-I))},$ $P=P_{N(W^{1/2}(AX_0-I))}$ and $Q=P_{N((W^{1/2}(AX_0-I))^{*})}.$

\vspace{0,3cm}
If $p=1,$ by Theorem \ref{TeoD} it holds, for every $\phi \in [0, 2\pi)$
$$0\leq D_{\phi}F_1(X_0,Y)=Re \ [e^{i\phi} tr(U^{*}W^{1/2}AY)] + \Vert QW^{1/2}AYP\Vert_{1},  \mbox{ for every } Y \in L(\HH).$$
Considering a suitabe $\phi$ for each $Y \in L(\HH),$ we get
$$\vert tr(U^{*}W^{1/2}AY) \vert \leq \Vert QW^{1/2}AYP\Vert_{1}, \mbox{ for every } Y \in L(\HH).$$
Observe that $R(Q)=N(U^{*})$ and $R(P)=N(U),$ therefore $U^{*}Q=PU^{*}=0.$ 

Let $Y \in L(\HH)$ then $\vert tr(U^{*}W^{1/2}AY) \vert  = \vert tr((I-P)U^* W^{1/2}AY) \vert = \vert tr(U^* W^{1/2}AY(I-P)) \vert \leq \Vert QW^{1/2}AY(I-P)P\Vert_{1}=0.$
Then
$$tr(U^{*}W^{1/2}AY)=0, \mbox{ for every } Y \in L(\HH).$$
Therefore 
$$U^{*}W^{1/2}A=0.$$
Hence
$$R(W^{1/2}A)\subseteq N(U^{*})=N((AX_0-I)^{*}W^{1/2}).$$ Therefore
$$(AX_0-I)^{*}W^{1/2}(W^{1/2}A)=0,$$ or equivalently
$$A^{*}WAX_0=A^{*}W, $$ and by Proposition \ref{Teo2} and Theorem \ref{Teo1}, the pair $(W,R(A))$ is compatible.  

\vspace{0,3cm}
If $1 < p < \infty,$ by Theorem \ref{TeoD} it holds, for every $\phi \in [0, 2\pi)$
$$0\leq D_{\phi}F_p(X_0,Y)=p \ Re \ [e^{i\phi} tr ( \vert W^{1/2}(AX_0-I) \vert^{p-1} U^{*} W^{1/2}AY)], \mbox{ for every }Y \in L(\HH).$$
Considering a suitable $\phi$ and $Y,$ it follows that
$$\vert W^{1/2}(AX_0-I) \vert^{p-1} U^{*} W^{1/2}A=0.$$ 
Since $N(\vert W^{1/2}(AX_0-I) \vert^{p-1})=N(\vert W^{1/2}(AX_0-I) \vert)$ it holds that
$$\vert W^{1/2}(AX_0-I) \vert U^{*} W^{1/2}A=0,$$ and therefore
$$A^{*}WAX_0=A^{*}W.$$ By Proposition \ref{Teo2} and Theorem \ref{Teo1}, the pair $(W,R(A))$ is compatible. 

\vspace{0,3cm}
Finally, if $X_0 \in L(\HH)$ minimizes problem \eqref{eq1}, we have proven that $X_0$ is a solution of the normal equation and by Proposition \ref{Teo2}, it is a $W$- inverse of $A$. Conversely, if $X_0$ is a $W$-inverse of $A,$ then by Proposition \ref{Teo2}, $X_0$ minimizes equation \eqref{eq2}, and by Proposition \ref{Prop Nashed}, it minimizes equation \eqref{eq1}.
\hfill \qed
\end{Dem}

\vspace{0,3cm}
\begin{obs} Let $\vert \vert \vert \cdot \vert \vert \vert$ be any unitarily invariant norm on a non-zero ideal $\mathcal{J}$ of $L(\HH).$ Given $W \in L(\HH)^{+}$ such that $W^{1/2} \in \mathcal{J},$ consider the norm associated to $W$ given by  $$\Vert X \Vert _{W}=\vert \vert \vert W^{1/2}X \vert \vert \vert,$$ for $X\in L(\HH).$
Let $A \in CR(\HH),$ if  the pair $(W, R(A))$ is compatible, then there exists a minimum of the set $\{ \Vert AX-I \Vert_{W}: X \in L(\HH) \}$  and 
$$\underset{X \in L(\HH)}{min} \Vert AX-I \Vert_{W}=\vert \vert \vert W_{/R(A)}^{1/2} \vert \vert \vert.$$
In particular if $\mathcal{J}=L(\HH)$ and we consider the operator norm $\Vert \cdot \Vert,$ the remark follows. \end{obs}

In fact, if the pair $(W, R(A))$ is compatible, by Theorem \ref{Teo1}, there exists $X_0\in L(\HH)$ such that $F(X_0)=min_{X\in L(\HH)}F(X)=W_{/R(A)},$ i.e. $$W_{/R(A)}=F(X_0)\leq F(X), \mbox{ for every } X\in L(\HH).$$ Since $W^{1/2} \in \mathcal{J},$  by Proposition \ref{Prop Nashed}, 
 $$\vert \vert \vert {W_{/R(A)}^{1/2}}\vert \vert \vert =\vert \vert \vert W^{1/2}(AX_0-I) \vert \vert \vert \leq \vert \vert \vert W^{1/2}(AX-I) \vert \vert \vert=\Vert AX-I \Vert_{W}, \mbox{ for every } X\in L(\HH),$$
then $$\underset{X \in L(\HH)}{min} \Vert AX-I \Vert_{W} = \Vert AX_0-I \Vert_{W} = \vert \vert \vert W_{/R(A)}^{1/2} \vert \vert \vert.$$   

\section{Weighted least squares problems II}

In this section we study the following problem:
given $A \in CR(\HH),$ $B \in L(\HH)$ and $W \in L(\HH)^{+}$ such that $W^{1/2} \in S_p$ for some $p$ with $1 \leq p < \infty,$ analize the existence of
\begin{equation}
\underset{X \in L(\HH)}{min} \Nphissw{AX-B}. \label{eq3}
\end{equation}

To study problem \eqref{eq3} we introduce the following associated problem: given $A \in CR(\HH),$  $B \in L(\HH),$ $W \in L(\HH)^{+}$ and $G:L(\HH) \ra L(\HH),$ $$G(X)=(AX-B)^{*}W(AX-B),$$
analize the existence of
\begin{equation}
\underset{X \in L(\HH)}{inf} G(X), \label{eq4}
\end{equation}

\begin{lema} \label{Prop31}Let $A \in CR(\HH),$ $B \in L(\HH)$ and  $W \in L(\HH)^{+},$ then the set $\{B^{*}E^{*}WEB: \ E^2=E, \ N(E)=R(A)\}$ has an infimum and
$$B^{*}W_{/R(A)}B=inf \{B^{*}E^{*}WEB: \ E^2=E, \ N(E)=R(A)\}.$$\end{lema}

\begin{Dem} If $W$ is invertible, then the pair $(W,R(A))$ is compatible and, by Theorem \ref{TeoShorted2},
$$W_{/ R(A)}=\mbox{min } \{ E^{*}WE: E^2=E, \ N(E)=R(A)\}.$$
In this case there is a unique projection $E_0$ where the minimum is attained (see \cite{CMSSzeged}), i.e. $W_{/ R(A)}=E_0^{*}WE_0.$
Then 
$$B^{*}W_{/R(A)}B=B^{*}E_0^{*}WE_0B\leq B^{*}E^{*}WEB,$$ for every projection $E$ with $N(E)=R(A),$ and then
$$min \{B^{*}E^{*}WEB: \ E^2=E, \ N(E)=R(A)\}=B^{*}W_{/R(A)}B.$$

For a non-invertible $W \in L(\HH)^{+}$, by Theorem \ref{TeoShorted}, it always hold that
$$B^{*}W_{/R(A)}B \leq B^{*}E^{*}WEB, \mbox{ for every projection } E \mbox{ such that } N(E)=R(A).$$
Therefore $B^{*}W_{/R(A)}B$ is a lower bound of $\{B^{*}E^{*}WEB: \ E^2=E, \ N(E)=R(A)\}.$

If $C \geq 0$ is any other lower bound for $\{B^{*}E^{*}WEB: \ E^2=E, \ N(E)=R(A)\},$ then for any $\varepsilon > 0,$ and any projection $E\in L(\HH)$ with $N(E)=R(A),$ we have $$C\leq B^{*}E^{*}WEB \leq B^{*}E^{*}(W+\varepsilon I) EB.$$ 
Since $W+\varepsilon I$ is invertible, it follows that $C\leq B^{*}(W+\varepsilon I)_{/R(A)}B,$ and since $\varepsilon$ is arbitrary, by \cite[Cor.~2]{Shorted2}, we conclude that $C\leq B^{*}W_{/R(A)}B.$ \hfill \qed
\end{Dem}

\vspace{0,3cm}
\begin{prop} \label{Teo11} Let $A \in CR(\HH),$ $B \in L(\HH)$ and  $W \in L(\HH)^{+},$ then the infimum of problem \ref{eq4} exists and
$$\underset{X \in L(\HH)}{inf} G(X)=B^{*}W_{/R(A)}B.$$\end{prop}

\begin{Dem} Following the same idea as in Proposition \ref{Teo7}, let $X\in L(\HH), $ then
$$(B-AX)^{*}W(B-AX)=(B-AX)^{*}W_{/R(A)}(B-AX)+(B-AX)^{*}W_{R(A)}(B-AX)=$$ $$=B^{*}W_{/R(A)}B+(B-AX)^{*}W_{R(A)}(B-AX)\geq B^{*}W_{/R(A)}B,$$
because $R(A) \subseteq N(W_{/R(A)}).$ Hence $B^{*}W_{/R(A)}B$ is a lower bound of $G(X).$ If $C\geq 0$ is any other lower bound of $G(X),$ then $$C\leq G(X), \mbox { for every } X \in L(\HH).$$
In particular, $$C\leq B^{*}E^{*}WEB,$$ where $E$ is any projection such that $N(E)=R(A);$
in fact $R((I-E)B) \subseteq R(I-E) = N(E) =R(A),$ then by Theorem \ref{Teo8}, there exists $X_0 \in L(\HH),$ such that $(I-E)B=AX_0,$ i.e. $(-EB)=AX_0-B.$

Therefore by Lemma \ref{Prop31}
$$C \leq inf \{B^{*}E^{*}WEB : \ E^{2}=E, \ N(E)=R(A) \}= B^{*}W_{/R(A)}B.$$ 
Thus, $$B^{*}W_{/R(A)}B=\underset{X \in L(\HH)}{inf} G(X).$$ 
\hfill \qed
\end{Dem}

\vspace{0,3cm}
\begin{thm}  \label{Teo5} Let $A \in CR(\HH),$ $W \in L(\HH)^{+}$ and $B\in L(\HH)$. Problem \ref{eq4} has a minimum, i.e., there exists $X_0\in L(\HH)$ such that $$\underset{X \in L(\HH)}{min}G(X)=G(X_0)=B^{*}W_{/R(A)}B$$ if and only if  
$R(B) \subseteq R(A) + R(A)^{\perp_{W}}.$ \end{thm}

\begin{Dem}
Suppose problem \ref{eq4} has a minimum and let $y \in R(B),$ then there exists $x\in\HH$ such that $y=Bx.$ If $X_0\in L(\HH)$ is such that $G(X_0)=min_{X\in L(\HH)} G(X)$ then $$\PI{G(X_0)x}{x}\leq \PI{G(X)x}{x}, \mbox{ for every } X\in L(\HH),$$ or equivalently, $$\NWW{(AX_0-B)x}\leq\NWW{(AX-B)x}, \mbox{ for every } X\in L(\HH).$$ Let $u_0=X_0x$ and let $z\in \HH$ be arbitrary, then there exists $X\in L(\HH)$ such that $z=Xx$. Therefore, $$\NWW{Au_0-Bx}\leq\NWW{Az-Bx}, \mbox{ for every } z\in \HH.$$  Therefore
$u_0$ is a $W-$LSS of $Az=Bx,$ then by Theorem \ref{thmWLSS}
$$y=Bx \in R(A) + R(A)^{\perp_{W}},$$ concluding that $R(B) \subseteq R(A) + R(A)^{\perp_{W}}.$

Conversely, if $R(B) \subseteq R(A) + R(A)^{\perp_{W}},$ by Theorem \ref{thmWinversa}, the operator $A$ admits a $W$-inverse in $R(B).$ Let $X_0$ be a $W$-inverse of $A$ in $R(B)$, then
$$\NWW{AX_0x-Bx} \leq \NWW{Az-Bx}, \mbox{ for every } x, z \in \HH.$$   
In particular, given $X \in L(\HH),$ consider  $z=Xx.$ Then for every $x\in\HH,$
$$\NWW{AX_0x-Bx}\leq \NWW{AXx-Bx}.$$ 
Hence, $$\NWW{AX_0x-Bx}\leq \NWW{AXx-Bx}, \mbox{for every } x\in \HH \mbox { and for every } X \in L(\HH),$$ or equivalently $$G(X_0)\leq G(X), \mbox{ for every } X\in L(\HH),$$ therefore the set $\{G(X) : X\in L(\HH)\}$ admits a minimum element.\hfill \qed
\end{Dem}

\vspace{0,3cm} If $R(B) \subseteq R(A) + R(A)^{\perp_{W}},$ then by Theorem \ref{Teo5}, problem \ref{eq4} attains a minimum, more precisely we proved that every $W$-inverse of $A$ in $R(B)$ minimizes $G(X),$ i.e. if $V_0 \in L(\HH)$ is a $W$-inverse of $A$ in $R(B),$ then $G(V_0)=B^{*}W_{/R(A)}B.$ Consider the set $$M_B=\{X\in L(\HH): G(X)=B^{*}W_{/R(A)}B\}.$$ The next proposition gives a characterization of the elements of $M_B.$ 

\begin{prop} \label{Teo6} Let $A \in CR(\HH),$ $B\in L(\HH)$ and $W \in L(\HH)^{+}.$ If $R(B) \subseteq R(A) + R(A)^{\perp_{W}}$ then the following conditions are equivalent:
\begin{itemize}
\item [i)]  $V \in M_B,$ i.e. $G(V)=\underset{X \in L(\HH)}{min} G(X),$
\item [ii)] $V$ is a $W$-inverse of $A$ in $R(B),$ 
\item [iii)] $V$ is a solution of the normal equation $$A^{*}W(AX-B)=0.$$
\end{itemize}\end{prop} 
\begin{Dem}
$i) \Leftrightarrow ii):$  It follows from the proof of Theorem \ref{Teo5}.

$ii) \Leftrightarrow iii):$ It was proven in Theorem \ref{thmWinversa}.\hfill \qed\end{Dem}

\vspace{0,3cm}

\begin{thm} \label {Teo31} Let $A \in CR(\HH),$ $B\in L(\HH)$ and $W \in L(\HH)^{+},$ such that $W^{1/2} \in S_p,$ for some $p$ with $1 \leq p < \infty$. Problem \eqref{eq3} admits a minimum if and only if $R(B) \subseteq R(A) + R(A)^{\perp_{W}}.$

In this case
$$\underset{X \in L(\HH)}{min}\Nphissw {AX-B} = \Nphiss {W_{/R(A)}^{1/2}B}.$$ 
Moreover, $X_0 \in L(\HH)$ satisfies
$$\Nphissw {AX_0-B} = \Nphiss {W_{/R(A)}^{1/2}B},$$
if and only if $X_0$ is a $W$-inverse of $A$ in $R(B).$
\end{thm}
\begin{Dem} If $R(B) \subseteq R(A) + R(A)^{\perp_{W}}$, by Theorem \ref{Teo5}, there exists $X_0\in L(\HH)$ such that $G(X_0)=\underset{X \in L(\HH)}{min}G(X)=B^{*}W_{/R(A)}B,$ i.e. $$G(X_0)=B^{*}W_{/R(A)}B\leq G(X), \mbox{ for every } X\in L(\HH).$$ 
Since $W^{1/2} \in S_p,$ by Proposition \ref{Prop Nashed} it holds that  $$\Nphiss {W_{/R(A)}^{1/2}B}=\Nphiss{W^{1/2}(AX_0-B)}=\Nphissw{AX_0-B}\leq \Nphissw{AX-B}, \mbox{ for every } X\in L(\HH),$$
then $$\underset{X \in L(\HH)}{min} \Nphissw {AX-B} = \Nphissw{AX_0-B}=\Nphiss {W_{/R(A)}^{1/2}B}.$$ 

The converse can be proven in a similar way as in Theorem \ref{Teo3}.

\vspace{0,3cm}
Finally, if $X_0 \in L(\HH)$ minimizes \eqref{eq3}, we have proven that $X_0$ is a solution of the normal equation and by Proposition \ref{Teo6}, it is a $W$- inverse of $A$ in $R(B)$. Conversely, if $X_0$ is a $W$-inverse of $A$ in $R(B),$ then by Proposition \ref{Teo6}, $X_0$ minimizes \ref{eq4}, and by Proposition \ref{Prop Nashed}, it minimizes \eqref{eq3}.
\hfill \qed
\end{Dem}

\vspace{0,3cm}
\begin{obs} $i)$ Let $A \in CR(\HH)$ and  $W \in L(\HH)^{+}$. Problem \ref{eq4} has a minimum for all $B \in L(\HH)$ if and only if the pair $(W, R(A))$ is compatible. If the pair $(W, R(A))$ is compatible then for every $B \in L(\HH),$ we have $R(B) \subseteq \HH=R(A)+R(A)^{\perp_W}$ (see Theorem \ref{Comp 1}), and by Theorem \ref{Teo5}, problem \ref{eq4} has a minimum for all $B \in L(\HH).$ The converse follows from Theorem \ref{Teo1} taking $B=I.$

\vspace{0,3cm}
$ii)$ Let $A \in CR(\HH),$ $B \in L(\HH)$ with $\overline{R(B))}=\HH$ and $W \in L(\HH)^{+}.$ If problem \ref{eq4} 
has a minimum, then the pair $(W, R(A))$ is cuasi-compatible, i.e., there exists a closed (densely defined) projection $Q$ with $R(Q)=R(A)$ and $W$- symmetric, i.e.,  $WQx = Q^{*} Wx,\mbox{ for every } x \in \mathcal{D}(Q),$ the domain of $Q,$ see \cite{Cuasi}.

In fact, if problem \ref{eq4} has a minimum, then by Theorem \ref{Teo5}, $R(B) \subseteq R(A)+R(A)^{\perp_W},$ therefore $\HH=\overline{R(B)} \subseteq \overline{R(A)+R(A)^{\perp_W}}.$ Let $\N=R(A)\cap R(A)^{\perp_W}.$ Note that $R(A)+R(A)^{\perp_W}=R(A)\dot{+} R(A)^{\perp_W} \cap \N^\perp$ and define $Q=P_{R(A) {\mathbin{\!/\mkern-3mu/\!}} R(A)^{\perp_W} \ominus \N}.$ Then $Q$ is a closed densely defined projection. By \cite[Prop~2.2]{Cuasi}, $Q$ is $W$- symmetric and the pair $(W, R(A))$ is cuasi-compatible.
\end{obs}

\section*{Acknowledgements}
The authors gratefully acknowledge the helpful suggestions made by the reviewer.

Maximiliano Contino was supported by Peruilh fundation and CONICET PIP 0168. Juan I. Giribet was partially supported by CONICET PIP 0168 and UBACyT2014.  A. Maestripieri was partially supported by CONICET PIP 0168.

\end{document}